\documentclass{amsart}

\usepackage[all]{xy}   
\usepackage{hyperref, graphicx} 
\usepackage[usenames,dvipsnames]{xcolor}
\usepackage{amssymb}

\theoremstyle{plain}
\newtheorem{lem}{Lemma}[section]
\newtheorem{cor}[lem]{Corollary}
\newtheorem{thm}[lem]{Theorem}

\theoremstyle{definition}
\newtheorem{ex}[lem]{Example}
\newtheorem{rem}[lem]{Remark}
\newtheorem{dfn}[lem]{Definition}

\newcommand{\Z}{\mathbb{Z}}               

\newcommand{\oedge}[1]{\ar@{.}[#1]}      
\newcommand{\edge}[1]{\ar@{-}[#1]}      
\newcommand{\medge}[1]{\ar@{=}[#1]|{>}}
\newcommand{\dedge}[1]{\ar[#1]}   
\newcommand{\lulab}[1]{\ar@{}[l]_<<{{}_{#1}}}   
\newcommand{\node}{*=0{\circ}}                       
\newcommand{\enode}{*=0{\circledcirc}}  

\begin{document}

\title[Canonical dimension and unimodular degree]{The canonical dimension of a semisimple group and the unimodular degree of a root system}

\author{Kirill Zainoulline}
\address[Kirill Zainoulline]{Department of Mathematics and Statistics, University of Ottawa, 150 Louis-Pasteur, Ottawa, ON, K1N 6N5, Canada}
\email{kirill@uottawa.ca}
\thanks{Partially supported by NSERC Discovery grant RGPIN-2015-04469}

\begin{abstract}
We produce a short and elementary algorithm to compute an upper bound for the canonical dimension of a spit semisimple linear algebraic group.
Using this algorithm we confirm previously known bounds by Karpenko and Devyatov 
as well as we produce new bounds (e.g. for groups of types $F_4$, adjoint $E_6$, for some semisimple groups).
\end{abstract}

\subjclass[2010]{14L17, 14C25}
\keywords{canonical dimension, algebraic group, divided-difference operator}

\maketitle

Let $G$ be a split semisimple linear algebraic group over a field $F$. 
The canonical dimension $cd(G)$ of $G$ was introduced by Berhuy-Reichstein in \cite{BR05}.
Since then it has turned into an important invariant of $G$-torsors which has numerous applications in the geometric theory of quadratic forms and twisted flag varieties (see \cite{Ka10}) 
as well as in the representation theory (see \cite{KM08}).
The canonical dimension was computed for most simple groups which possess a unique torsion prime $p$ by Karpenko in~\cite{Ka08, Ka06, Ka05}. 
Its $p$-local version $cd_p(G)$ (which also gives a natural lower bound for $cd(G)$) was computed for all simple groups by Karpenko-Merkurjev in~\cite{KM06} and by the author in~\cite{Za07}.
If $G$ has several torsion primes the only known case of $PGL_6$ was treated in \cite{CKM}. 

In the present notes we provide a short and elementary algorithm to compute an upper bound for $cd(G)$ for all split semisimple groups. 
As an application, we obtain upper bounds computed by Karpenko using the Pieri formula for the intersection of effective divisors.
We also simplify Devyatov's computations in \cite{De17}, where such bounds were obtained for all simply-laced simply-connected simple groups using advanced combinatorial techniques.
In addition, we get new bounds  (e.g. for $F_4$, adjoint $E_6$, and some semisimple groups which are not direct-products of simple groups).
Instead of the Pieri-type formula our arguments rely on the Demazure formula~\cite[(9)]{De73} for the characteristic map 
and elementary properties of the divided-difference operators (see \cite[\S4]{De73}).

The paper is organized as follows: In section~\ref{sec:bcan} following \cite{Ka06,KM06} we recall how to obtain an upper bound for $cd(G)$ using the intersection of effective divisors in the complete flag variety.
In section~\ref{sec:unim} we introduce the notion of a unimodular degree of a root system and produce an elementary algorithm to estimate it.
Our main technical tools here are Theorem~\ref{thm:main} and Example~\ref{ex:c3}.
In the last section, using the unimodular degree, we provide various examples of computations. 
We summarize upper bounds for simply-connected simple groups in~\ref{subsec:sc}. In~\ref{subsec:nonsc} we explain how to obtain such bounds for non-simply-connected groups. 

\section{An upper bound for the canonical dimension}\label{sec:bcan}

We recall Karpenko's method (see \cite{Ka06} for details) to obtain an upper bound for the canonical dimension of $G$. 

\subsection{\it The canonical dimension of a variety.}\label{sub:cdvar}

Recall that
given a smooth projective variety $X$ over $F$, its canonical dimension $cd(X)$
is defined to be the minimum of $\dim Y$, 
where $Y$ runs over closed irreducible subvarieties of $X$ such that the base change $Y_{F(X)}$ has a rational point (see e.g. \cite[Definition~1.1]{Ka10}). 

Suppose we are given a non-negative cycle $\gamma=\sum_i a_i[Y_i] \in CH_N(X)$ of dimension $N$, where each $Y_i$ is a closed irreducible subvariety in $X$ and $a_i\ge 0$.
The base change to the function field $K=F(X)$ gives a non-negative cycle in $CH_N(X_K)$ that is
$\gamma_K=\sum_{i,j}a_i b_{i,j} [Z_{i,j}]$, where $[(Y_i)_K]=\sum_j b_{i,j} [Z_{i,j}]$, $b_{i,j}\ge 0$ and $Z_{i,j}$ are closed irreducible subvarieties in $X_K$.

Assume that the tangent bundle of $X_K$ is generated by global sections. Then any product of non-negative cycles in $CH(X_K)$ 
is represented by a non-negative cycle. So, if there exists an irreducible subvariety $Z$ in $X_K$ such that $\gamma_K\cdot [Z]=[pt_K]$ gives the class of a $K$-point, 
then $[Z_{i,j}]\cdot [Z]=m_{i,j} [pt_K]$, where $m_{i,j}\ge 0$ and $\sum a_ib_{i,j}m_{i,j}=1$. Because of the non-negativity, there are unique $i$ and $j$ with $a_i=b_{i,j}=m_{i,j}=1$ in which case
$Z_{i,j}\cap Z$ and, hence, $(Y_i)_K$ contains a $K$-point. By definition of the canonical dimension it then gives $cd(X)\le N$.

Summarizing the arguments above, to find the lower bound for $cd(X)$ we need
to find a non-negative cycle $\gamma$ on $X$ of minimal possible dimension such that its base change $\gamma_K$ has a `Poincar\'e dual' in $X_K$.

\subsection{\it The canonical dimension of a group.}\label{sub:cdgr}

We fix a split maximal torus $T$ of rank $n$ and a Borel subgroup $B$ of $G$ containing $T$. 
Let $E$ be a $G$-torsor over $F$ and let $E/B$ denote the respective twisted form of the variety of complete flags $G/B$.
By \cite[\S6, p.424]{KM06} the canonical dimension of the group $G$ is the supremum
$cd(G)=sup_E\{cd(E/B)\}$ taken over all $G$-torsors $E$ over all field extensions of $F$.
So our task reduces to constructing such non-negative  cycle $\gamma$ for all $E/B$'s.

Let $\Sigma$ be the semisimple root system associated to $T\subset B$ with positive roots $\Sigma^+$
and simple roots $\Pi$. Let $T^*$ be the group of characters.
We have $\Lambda_r\subset T^*\subset \Lambda_w$, where $\Lambda_r$ is the root lattice and $\Lambda_w$ is the weight lattice of $\Sigma$. 
Consider the characteristic map
\[c\colon Sym_{\Z}(T^*) \to CH(G/B_K)\quad \text{ given by } \lambda\mapsto c_1(L_\lambda),\]
where $\lambda\colon B\to Aut(V)=\mathbb{G}_m$ is a character extended to $B$ and $L_\lambda= G\times^B \mathbb{G}_m=G\times V/(gb^{-1},\lambda(b)v)$ is the associated line bundle. If $\lambda$ is dominant, i.e., $\lambda\in \Lambda_w^+\cap T^*$, where $\Lambda_w^+$ denote the set of dominant weights,
then $L_\lambda$ has a non-zero global section, and, hence, $c_1(L_\lambda)\in CH(G/B_K)$ is non-negative.

Since the tangent bundle of $E/B$ (as well as of $G/B_K$) is generated by global sections, any product $\gamma$ of non-negative cycles of codimension one in $E/B$ is also non-negative.
The restriction map on the Picard groups \[res\colon CH^1(E/B)\to CH^1(G/B_K)\] is injective and by \cite[Thm.6.4.(1)]{KM06} its image contains the image of the characteristic map.
So
$res(\gamma_\lambda)=c(\lambda)$, where $\gamma_\lambda \in CH^1(E/B)$ is non-negative if so is $c(\lambda)$.
Let $Sym^+_{\Z}(T^*)$ denote the non-negative linear span of products  of $\lambda\in \Lambda_w^+\cap T^*$.
Then any homogeneous polynomial $p \in Sym^+_{\Z}(T^*)$ of degree $N$ gives rise to a non-negative cycle $\gamma\in CH_N(E/B)$ such that $\gamma_K=c(p)$.

Consider an additive basis of $CH(G/B_K)$ given by the classes of the Schubert varieties (closures of Bruhat cells) $V^w$, where $w\in W$ runs through all elements of the Weyl group $W$. 
So we can write $\gamma_K=\sum_{w\in W} a_wV^w$, where $a_w\ge 0$.
Suppose now that there is $w$ such that $a_w=1$ (we call such $\gamma_K$ unimodular). Then taking $Z:=V^{w_0w}$ to be the Poincar\`e dual of $V^w$ we obtain the desired condition
$\gamma_K\cdot [Z]=[pt_K]$. Since $G/B$ is a cellular space over $F$, we may identify $CH(G/B_K)$ with $CH(G/B)$. Hence, $cd(G)\le N$, where $N$ is the smallest dimension of a non-negative unimodular cycle in $CH(G/B)$.

\section{Unimodular degree of a root system}\label{sec:unim}

In the present section we introduce the (purely combinatorial) concept of a unimodular degree $ud(\Sigma)$ of a root system. 
We produce an elementary algorithm to obtain lower bounds for $ud(\Sigma)$ for all semisimple finite root systems. 
Our enumeration of roots follows Bourbaki~\cite{BRB}.

\subsection{\it The unimodular degree.} Let $\Sigma$ be a semisimple finite root system of rank $n$ with a subset of simple roots $\{\alpha_1,\ldots,\alpha_n\}$.
Let $C=(c_{i,j})_{i,j}$ be the Cartan matrix of $\Sigma$. By definition (following \cite{BRB}) we have \[c_{i,j}=\alpha_j^\vee(\alpha_i)=\tfrac{2\alpha_i\cdot \alpha_j}{\alpha_j\cdot \alpha_j},\;\text{ where }c_{i,j}\le 0\text{ if }i\neq j.\]
Moreover, $c_{i,j}\neq 0$, $i\neq j$ if and only if the $i$th and the $j$th vertices of the Dynkin diagram
are connected by an edge. 

Let $x_1,\ldots,x_n$ be the fundamental weights of $\Sigma$. By definition we have $\alpha_i=\sum_{j=1}^nc_{i,j}x_j$.
Let $S=\Z[x_1,\ldots,x_n]$ denote the polynomial ring in $x_i$s. 
Set
\[
y_i:=x_i-\alpha_i=-c_{i,1}x_1\ldots-c_{i,i-1}x_{i-1}-x_i-c_{i,i+1}x_{i+1}-\ldots-c_{i,n,}x_n,
\]

Consider the divided-difference $\Z$-linear operator (cf.~\cite[\S4]{De73})
\[
\partial_i(p)=(p-s_i(p))/(x_i-y_i),\quad p\in S,
\]
where $s_i$ is the $i$th simple reflection defined by
$
s_i(x_i)=y_i$ and $s_i(x_j)=x_j$ if $i\neq j$.
It maps a homogeneous polynomial of degree $d$ into a homogeneous polynomial of degree $d-1$. 
By definition, we have
$
\partial_i (x_j)=\delta_{i,j}^{Kr},
$
and $\partial_i$ satisfies the twisted Leibniz rule:
\[
\partial_i(pq)=\partial_i(p)q+s_i(p)\partial_i(q),\; p,q\in S.
\]

Given a reduced word $w$ of length $l=\ell(w)$ in the Weyl group $W$ of $\Sigma$, i.e., $w=s_{i_1}\ldots s_{i_l}$ we set 
\[\partial_w:=\partial_{i_1}\ldots \partial_{i_l} \text{ (with respect to the usual composition of operators)}.
\]
The composite $\partial_w$ does not depend on a choice of reduced word so it is well-defined. The operators $\partial_i$s
generate the nil-Coxeter ring of $\Sigma$. In particular, they satisfy the standard nil-Coxeter relations, e.g., \[\partial_i^2=0,\; (\partial_i\partial_j)^{m_{ij}}=1,\quad 1\le i<j\le n\] where $m_{ij}$ are the Coxeter exponents of $\Sigma$. 

\begin{dfn}
We define the unimodular degree of $\Sigma$ to be
\[
ud(\Sigma):=\max\,\{\deg p \mid p\text{ is a monic monomial},\; \exists w\in W: \partial_w(p)=1\}
\]
Here the maximum is taken over degrees of monic monomials $p$ such that $\partial_w(p)=1$ for some $w\in W$. Observe that $\partial_w(p)=1$ $\Longrightarrow$ $\ell(w)=\deg p$.
\end{dfn}

\subsection{\it Properties of divided-difference operators.}

We start with the following basic property:

\begin{lem}\label{lem:prod}
Given two homogeneous $p,q\in S$, if $p$ does not contain a monomial divisible by $x_i$, then $\partial_i(pq)=p\partial_i(q)$.
\end{lem}

\begin{proof} By induction on the degree of $p$ and the twisted Leibniz rule.
\end{proof}

\begin{cor}
If $i\notin \{i_1,\ldots,i_k\}$, then $\partial_i(x_{i_1}^{e_1}\ldots x_{i_k}^{e_k})=0$, $e_1,\ldots,e_k\ge 0$.
\end{cor}

\begin{cor}
Let $w=s_{j_1}\ldots s_{j_k}$, where $\{j_1,\ldots,j_k\}=\{i_1,\ldots,i_k\}$. 

Then $\partial_w(x_{i_1}\ldots x_{i_k})=1$.
In particular, $ud(\Sigma)\ge n$.
\end{cor}

We now look at powers of fundamental weights:

\begin{lem}\label{lem:power}
We have
$
\partial_i(x_i^e)=\sum_{m=1}^{e} x_i^{m-1} y_i^{e-m}$, $e\ge 1$.

In particular, $\partial_i(x_i^2)=x_i+y_i$ and
$\partial_i(x_i^3)=x_i^2+x_iy_i +y_i^2$.
\end{lem}

\begin{proof} By induction on $e$ using the Leibniz rule
\[
\partial_i(x_i^e)=\partial_i(x_i)x_i^{e-1}+s_i(x_i)\partial_i(x_i^{e-1})
=x_i^{e-1}+y_i\partial_i(x_i^{e-1}).\qedhere
\]
\end{proof}

\begin{ex}\label{ex:c3}
By direct computations for the root system of type $C_3$ we obtain
\[
\partial_3(x_3)=\partial_2\partial_3\partial_2(x_2^3)=\partial_1\partial_2\partial_3\partial_2\partial_1(x_1^5)=1.
\]
So for the monomial $p=x_1^5x_2^3x_3$ of degree $9$ we get
\[\partial_1\partial_2\partial_3\partial_2\partial_1(\partial_2\partial_3\partial_2(\partial_3(x_1^5x_2^3x_3)))=1\]
which implies $ud(C_3)\ge 9$. Since the element $s_1s_2s_3s_2s_1s_2s_3s_2s_3$ is of maximal length, we obtain that $ud(C_3)=9$.

Observe that by the same reasons $ud(C_2)=4$ which corresponds to $p=x_2^3x_3$.
\end{ex}

\begin{dfn}
We say that a $k$-tuple $(i_1,\ldots,i_k)$, $k\ge 2$ of distinct indices is a $1$-chain if \[c_{i_k,i_{k-1}}=\ldots =c_{i_2,i_1}=-1.\]
We set $(i_1)$ to be a $1$-chain by default.
\end{dfn}
Observe that a $1$-chain corresponds to a subchain of the Dynkin diagram of $\Sigma$ so that the length of $\alpha_{i_1}$ is greater or equal to the length of $\alpha_{i_k}$.

\begin{thm}\label{thm:main}
Let $(i_1,\ldots,i_k)$ be a $1$-chain, and let $p$ be a monomial coprime to $x_{i_1}\ldots x_{i_k}$. 
Then 
\begin{itemize}
\item[(i)] $\partial_{i_1}\ldots \partial_{i_k}(x_{i_k}^e p)=0$ if $0\le e<k$; and \item[(ii)] $\partial_{i_1}\ldots \partial_{i_k}(x_{i_k}^k p)=p$.
\end{itemize}
\end{thm}

\begin{proof} We prove both (i) and (ii) by induction on $k$. The case $k=1$ follows by Lemma~\ref{lem:prod} and its corollaries. 

(i)  
Assume $k\ge 2$.
By the Leibniz rule and Lemma~\ref{lem:power}
\begin{align*}
\partial_{i_1}\ldots \partial_{i_k}(x_{i_k}^e p)&=\partial_{i_1}\ldots \partial_{i_{k-1}}\big((\sum_{m=1}^e x_{i_k}^{m-1}y_{i_k}^{e-m})p\big)\\
&=\sum_{m=1}^{e}\partial_{i_1}\ldots \partial_{i_{k-1}} \big(y_{i_k}^{e-m} (x_{i_k}^{m-1}p)\big),
\end{align*}
where $y_{i_k}$ contains a summand $-c_{i_{k},i_{k-1}}x_{i_{k-1}}=x_{i_{k-1}}$ but it does not contain summands $x_{i_1},\ldots,x_{i_{k-2}}$ (being a $1$-chain implies that $c_{i_k,i_1}=\ldots =c_{i_{k},i_{k-2}}=0$).
After collecting the coefficients at the powers of $x_{i_{k-1}}$ in $y_{i_k}^{e-m} (x_{i_k}^{m-1}p)$ we finish the proof of (i) by induction.

(ii) By the Leibniz rule and Lemma~\ref{lem:power} we obtain
\[
\partial_{i_1}\ldots \partial_{i_k}(x_{i_k}^k p)=\partial_{i_1}\ldots \partial_{i_{k-1}}\big((y_{i_k}^{k-1}+y_{i_k}^{k-2}x_{i_k}+\ldots + x_{i_k}^{k-1})p\big).
\]
Since $y_{i_k}=-c_{i_{k},i_{k-1}}x_{i_{k-1}}+q_k=x_{i_{k-1}}+q_k$, where $q_k$ does not contain any monomial with $x_{i_j}$, $j=1\ldots k-2$, collecting the coefficients at the powers of $x_{i_{k-1}}$ we can rewrite the latter expression as
\[
\partial_{i_1}\ldots \partial_{i_{k-1}}\big((x_{i_{k-1}}^{k-1}+\sum_{m=1}^{k-1}x_{i_{k-1}}^{k-1-m} f_m)p\big),
\]
where each $f_m$ does not contain monomials with $x_{i_j}$, $j=1\ldots k-1$.
The proof of (ii) is then finished by induction and by (i).
\end{proof}

\subsection{\it The lower bounds and 1-chains.}
Our goal now is to find a lower bound for $ud(\Sigma)$ using 1-chains of Theorem~\ref{thm:main}. The following example demonstrates how to choose `the best possible' sequence of 1-chains.

\begin{ex} Consider the root system of type $E_8$. The following sequence of 1-chains produces a monomial $p$ of maximal possible degree $34$.

\begin{center}
{\scriptsize
\begin{tabular}{rcl}
$\xymatrix@C=1.5em@R=1ex{
 &\node\lulab{1}\oedge{r}&\node\lulab{3}\oedge{r}& \node \lulab{4} \dedge{d} & \node \lulab{5} \dedge{l} &
\node \lulab{6} \dedge{l} & \node \lulab{7}
\dedge{l} & \node \lulab{8} \dedge{l} \\
&&& \enode \lulab{2}  &&&&
}$ & $\quad\mapsto\quad$ & $\partial_8\partial_7\partial_6\partial_5\partial_4\partial_2(x_2^6)=1$, $p:=x_2^6$ \\
$\xymatrix@C=1.5em@R=1ex{
 &\enode\lulab{1} &\node\lulab{3}\dedge{l}& \node \lulab{4} \dedge{l} & \node \lulab{5} \dedge{l} &
\node \lulab{6} \dedge{l} & \node \lulab{7}
\dedge{l} & \node \lulab{8} \dedge{l} \\
&&& \node \lulab{2} \oedge{u} &&&&
}$ & $\quad\mapsto\quad$ & $\partial_8\partial_7\partial_6\partial_5\partial_4\partial_3\partial_1(x_1^7p)=p$, $p:=x_1^7x_2^6$ \\
$\xymatrix@C=1.5em@R=1ex{
 &\node\lulab{1} &\enode\lulab{3}\oedge{l}& \node \lulab{4} \dedge{l} & \node \lulab{5} \dedge{l} &
\node \lulab{6} \dedge{l} & \node \lulab{7}
\dedge{l} & \node \lulab{8} \dedge{l} \\
&&& \node \lulab{2} \oedge{u} &&&&
}$ & $\quad\mapsto\quad$ & $\partial_8\partial_7\partial_6\partial_5\partial_4\partial_3(x_3^6p)=p$, $p:=x_3^6x_1^7x_2^6$ \\
 & \vdots & \\
$\xymatrix@C=1.5em@R=1ex{
 &\node\lulab{1} &\node\lulab{3}\oedge{l}& \node \lulab{4} \oedge{l} & \node \lulab{5} \oedge{l} &
\node \lulab{6} \oedge{l} & \enode \lulab{7}
\oedge{l} & \node \lulab{8} \dedge{l} \\
&&& \node \lulab{2} \oedge{u} &&&&
}$ & $\quad\mapsto\quad$ & $\partial_8\partial_7(x_7^2p)=p$, $p:=x_7^2x_6^3x_5^4x_4^5x_3^6x_1^7x_2^6$ \\
$\xymatrix@C=1.5em@R=1ex{
 &\node\lulab{1} &\node\lulab{3}\oedge{l}& \node \lulab{4} \oedge{l} & \node \lulab{5} \oedge{l} &
\node \lulab{6} \oedge{l} & \node \lulab{7}
\oedge{l} & \enode \lulab{8} \oedge{l} \\
&&& \node \lulab{2} \oedge{u} &&&&
}$ & $\quad\mapsto\quad$ & $\partial_8(x_8p)=p$, $p:=x_8x_7^2x_6^3x_5^4x_4^5x_3^6x_1^7x_2^6$ \\
\end{tabular}
}
\end{center}

\ 

Here we start with $p:=1$ and on each step multiply it by the monomial $x_i^k$, where $k$ is the number of vertices in the respective 1-chain ending at the $i$th vertex (circled).
\end{ex}

\begin{ex}
Similarly to the previous example, by looking only at 1-chains we obtain the following list of polynomials and the respective lower bounds (their degrees) for all simple root systems:

\begin{center}
{\scriptsize
\begin{tabular}{c|c|c|c}
Type & $p=$ &  $ud(\Sigma)\ge$ & $|\Sigma^+|$  \\
\hline
$A_n$ & $x_1^n x_2^{n-1}\ldots x_2^2x_1$ &  ${n+1 \choose 2}$ & ${n+1 \choose 2}$ \\
$B_n$ & $x_1x_2^2\ldots x_{n-1}^{n-1}x_n^n$  & ${n+1 \choose 2}$ & $n^2$\\
$C_n$ & $x_1^nx_2^{n-1}\ldots x_{n-1}^2x_n$  & ${n+1 \choose 2}$ & $n^2$\\
$D_n$ & $x_1x_2^2\ldots x_{n-2}^{n-2}x_{n-1}^{n-1}x_n^{n-1}$  & ${n+1 \choose 2}-1$ & $2{n \choose 2}$\\
$E_n$ & $x_1^{n-1}x_2^{n-2}x_3^{n-2}\ldots x_{n-1}^2 x_n$  & ${n+1 \choose 2}-2$ & 36, 63, 120 \\
$F_4$ & $x_1x_2^2x_3^3x_4^4$  & 10 & 24\\
$G_2$ & $x_1x_2^2$  & 3 & 6\\
\end{tabular}
}
\end{center}
\end{ex}

\begin{rem} Observe that the bound for type $A$ is sharp 
as it coincides with the maximal length of an element in $W$ that is $|\Sigma^+|$.
By the main result of \cite{De17} the  above bounds are also sharp for types $D$ and $E$.
\end{rem}

\begin{ex} Combining the equalities of Example~\ref{ex:c3} with 1-chains for the root system of type $F_4$ as follows:

\begin{center}
{\scriptsize
\begin{tabular}{rcl}
$\xymatrix@C=1.5em@R=1ex{
 &\node\lulab{1}\dedge{r}&\node\lulab{2}\medge{r}& \node \lulab{3} \dedge{r}  & \enode \lulab{4} 
}$ & 
$\quad\mapsto\quad$ & $\partial_4\partial_3\partial_2\partial_3\partial_4(x_4^5)=1$ (Example~\ref{ex:c3}), $p:=x_4^5$  \\
$\xymatrix@C=1.5em@R=1ex{
 &\node\lulab{1}\dedge{r}&\node\lulab{2}\medge{r}& \enode \lulab{3} \oedge{r}  & \node \lulab{4} 
}$ & 
$\quad\mapsto\quad$ & $\partial_3\partial_2\partial_3(x_3^3p)=1$ (Example~\ref{ex:c3}), $p:=x_3^3x_4^5$  \\
$\xymatrix@C=1.5em@R=1ex{
 &\node\lulab{1}\dedge{r}&\enode\lulab{2}\medge{r}& \node \lulab{3} \oedge{r}  & \node \lulab{4} 
}$ & 
$\quad\mapsto\quad$ & $\partial_1\partial_2(x_2^2)=1$ (1-chain), $p:=x_2^2x_3^3x_4^5$  \\
$\xymatrix@C=1.5em@R=1ex{
 &\enode\lulab{1}\oedge{r}&\node\lulab{2}\medge{r}& \node \lulab{3} \oedge{r}  & \node \lulab{4} 
}$ & 
$\quad\mapsto\quad$ & $\partial_1(x_1)=1$ (1-chain), $p:=x_1x_2^2x_3^3x_4^5$  \\
\end{tabular}
}
\end{center}
we improve the bound  as $ud(F_4)\ge \deg p=11$.

\end{ex}

\section{The unimodular degree and the canonical dimension}
In this section we relate $ud(\Sigma)$ with the canonical dimension $cd(G)$ for all semisimple groups $G$. We produce new bounds for type $F_4$ and adjoint groups of types $E_6$ 
as well as for some semisimple groups which are not direct-products of simple groups.

According to~\cite[(9)]{De73} we have the following formula for the characteristic map
\[
c(p)=\sum_{w\in W} \varepsilon \partial_w(p)[V^w],
\]
where in the notation of \cite[\S4]{De73} $H=CH(G/B)$, $S(M)=S$, $z_w=[V^w]$, $\varepsilon\colon S\to \Z$ is the augmentation map $x_i\mapsto 0$ and $D_w=\partial_w$.
By the very definition, $ud(\Sigma)$ is then the largest codimension of a unimodular cycle $c(p)$, where $p\in Sym^+_{\Z}(\Lambda_w)$. By section~\ref{sec:bcan} it follows that \[cd(G)\le dim(G/B)-ud(\Sigma)=|\Sigma^+|-ud(\Sigma).\]

\subsection{\it Simply-connected groups.}\label{subsec:sc}
If $G=G^{sc}$ is simply-connected (so $T^*=\Lambda_w$), then the lower bounds for $ud(\Sigma)$  
give upper bounds for $cd(G^{sc})$ which can be summarized as follows:

\begin{center}
{\scriptsize
\begin{tabular}{c|c|c|c|c|c|c|c}
 &  $B_n$ & $D_n$ & $E_6$ & $E_7$ & $E_8$ & $F_4$ & $G_2$ \\
  \hline
$cd(G^{sc})\le $  & ${n \choose 2}$  & ${n-1 \choose 2}$ & 17 & 37 & 86 & 13 & 3
\end{tabular}
}
\end{center}

\

Observe that $cd(G^{sc})=0$ for types $A$ and $C$, as in this case the torsion index of $G^{sc}$ is $1$. 

\begin{rem}
For the simply-connected group $G^{sc}$ of type $C_n$ using 1-chains only one obtains the bound $cd(G^{sc})\le {n \choose 2}$ (as in the $B_n$ case).
However, $cd(G^{sc})=0$ would follow if  $\partial_i\ldots \partial_{n-1}\partial_n\partial_{n-1}\ldots \partial_i(x_i^{2n-2i+1})=1$, $i=1\ldots n$ (similar to Example~\ref{ex:c3}).
\end{rem}

\begin{rem}
By \cite{Ka06} the above bounds are sharp for types $B_n$ and $D_{n+1}$ (which correspond to $Spin_{2n+1}$ and $Spin_{2n+2}$ groups)  if $n+1$ is a power of $2$, and
it is also sharp for $G_2$.
\end{rem}

\subsection{\it Non-simply-connected groups}\label{subsec:nonsc}
Let $G$ be a non-simply-connected group. Then $G=G^{sc}/Z$, where $Z$ is a non trivial central subgroup of $G^{sc}$.
The sublattice $T^*$ of the weight lattice $\Lambda_w$ can be identified with the kernel of the natural projection on characters $\pi\colon \Lambda_w \to Z^*$.
Observe that $Z^*$ is a finite abelian group which is cyclic for all simple groups except for adjoint groups of type $D$ of even ranks 

The upper bound for $cd(G)$ can be obtained by adjusting the 1-chain method (used in the simply-connected case) as follows.

We choose fundamental weights $x_{k_1},\ldots, x_{k_r}$ such that $\pi(x_{k_1}), \ldots,\pi(x_{k_r})$ generate $Z^*$,
and we replace each remaining fundamental weight $x_i$, $i\notin \{k_1,\ldots,k_r\}$ 
by the dominant weight 
\[
z_i=x_i+\sum_{j=1}^r a_{i,j}x_{k_j}\in T^*,\quad a_{i,j}\ge 0.\] 

Then we remove all $k_1,\ldots, k_r$-th vertices from the (semisimple) Dynkin diagram of $G$, 
and apply the 1-chain method as usual to the monomials in $z_i$s and the remaining Dynkin diagram.

\begin{ex} 
Consider the half-spin group $G=HSpin_{2n}$ that is of type $D_{n}$, $2\mid n$. The projection $\pi\colon \Lambda_w\to \Z/2\Z$ is then given by $x_i\mapsto i\mod 2$.
Removing $x_{n-1}$ (the $(n-1)$th vertex) and applying the 1-chain method to the remaining diagram of type $A_{n-1}$
gives 
$cd(G)\le |\Sigma^+|-ud(A_{n-1})\le {n \choose 2}$.
Observe that by \cite[Cor.~1.5]{Ka06} this bound is sharp if $n$ is a power of $2$.

For the adjoint group $G=PGO_{2n}$, $2\mid n$ the projection $\pi\colon \Lambda_w \to \Z/2\Z\times \Z/2\Z$ is given by $x_i\mapsto (i\mod 2,i\mod 2)$ for $i<n-1$, and $x_{n-1}\mapsto (1,0)$, $x_n\mapsto (0,1)$.
Removing both $x_{n-1}$ and $x_{n}$ we obtain the diagram of type $A_{n-2}$ which gives
$cd(G)\le |\Sigma^+|-ud(A_{n-2})\le \tfrac{(n-1)(n+2)}{2}$. Again, this bound becomes sharp if $n$ is a power of 2 by \cite[\S8.4]{KM06}. 
\end{ex}

\begin{ex}
Consider the adjoint group $G$ of type $E_6$. From \cite[Ch.VI \S4.12]{BRB} we see that $\pi\colon \Lambda_w \to \Z/3\Z$ maps
$x_1,x_5\mapsto 1$, $x_3,x_6\mapsto 2$ and $x_2,x_4\mapsto 0$. Removing $x_1$ we obtain the diagram of type $D_5$ which gives
$cd(G)\le |\Sigma^+|-ud(D_5)\le 22$.

Let $G$ be the adjoint group of type $E_7$. By \cite[Ch.VI \S4.11]{BRB} the projection $\pi\colon \Lambda_w \to \Z/2\Z$ is defined by
$x_1,x_3,x_4,x_6\mapsto 0$ and $x_2,x_5,x_7\mapsto 1$. Removing $x_2$ we obtain
$cd(G)\le |\Sigma^+|-ud(A_6)\le 42$. Observe that $42$ is also the dimension of the projective homogeneous variety $G/P_2$ corresponding to the generically split vertex on the Tits diagram in~\cite{Ti66}.
\end{ex}

\begin{ex} 
Let $G$ be a simply-connected simple group with the centre $Z$ corresponding to the root system $\Sigma$. 
Consider the direct product $\tilde G=G^m/Z$ of $m$ copies of $G$ modulo the diagonal subgroup $Z$. 
Then
\[
cd(\tilde G)\le m|\Sigma^+| - (m-1)ud(\Sigma)-ud(\Sigma'),
\]
where $\Sigma'$ is the root subsystem obtained by removing the respective simple roots $\alpha_{k_1},\ldots,\alpha_{k_r}$.
For instance, 
\[cd (SL_{n+1}^m/\mu_{n+1})\le n\quad \text{ and }\quad cd (E_6^2/\mu_3) \le 39.\]
\end{ex}

\bibliographystyle{alpha}

\end{document}